# Decentralized Demand Side Management with Rooftop PV in Residential Distribution Network

Hamidreza Sadeghian, Zhifang Wang
Department of Electric and Computer Engineering
Virginia Commonwealth University
Richmond, VA, USA
Email: {sadeghianh, zfwang} @ vcu.edu

*Abstract*—In the past extensive researches have been conducted on demand side management (DSM) program which aims at reducing peak loads and saving electricity cost. In this paper, we propose a framework to study decentralized household demand side management in a residential distribution network which consists of multiple smart homes with schedulable electrical appliances and some rooftop photovoltaic generation units. Each smart home makes individual appliance scheduling to optimize the electric energy cost according to the day-ahead forecast of electricity prices and its willingness for convenience sacrifice. Using the developed simulation model, we examine the performance of decentralized household DSM and study their impacts on the distribution network operation and renewable integration, in terms of utilization efficiency of rooftop PV generation, overall voltage deviation, real power loss, and possible reverse power flows.

*Index Terms*-- Demand Side Management, Distributed Generation, Rooftop Photovoltaic, Penalty factor, Optimization algorithm.

## NOMENCLATURE

| | |
|---|---|
| $u_a$ | Binary status of appliance $a$; 0 = *off*, 1 = *on* |
| $r_a$ | Rated power in kW for appliance $a$ |
| $A$ | Total number of appliances |
| $A_{int}$ | Total number of interruptible appliances |
| $T$ | Total number of time slots, $T=48$ |
| $t$ | Index for time slots |
| $D_a$ | During of operation time for appliance $a$ |
| $[s_a, f_a]$ | Allowable operational time range for appliance $a$ |
| $t_a^{st}$ | Operation start time of appliance $a$ |
| $t_a^{en}$ | Operation end time of appliance $a$ |
| $t_a^{st_{new}}$ | Operation start time of appliance $a$ after DSM shifting |
| $t_a^{st_{old}}$ | Operation start time of appliance $a$ before shifting |
| $P_{pv}(t)$ | PV output in kW at time slot $t$ |
| $\pi_e(t)$ | Electricity price (¢/kWh) at time slot $t$ |
| $\pi_p$ | Penalty price (¢/kWh) |
| $\Delta T_a$ | Time shifting for appliance $a$ |
| $MD$ | Maximum demand of household power consumption |
| $B$ | Total number of buses |
| $b$ | Index for buses |
| $V(t)$ | Vector of bus voltage at time slot $t$ |
| $P(t)$ | Vector of active load power at t |
| $Q(t)$ | Vector reactive power at t |
| $f_{AC}(.)$ | AC power flow equation |
| $P_{loss}(t)$ | Active power loss in distribution feeder |
| $C_e$ | Total electricity cost (¢/day) |
| $C_p$ | Total penalty cost (¢/day) |

## I. INTRODUCTION

With recent technology advances, solar photovoltaic has become one of the fastest-growing renewable energy sources in the U.S. [1]. However, high penetration of PV systems into the distribution networks may arise undesirable issues such as voltage fluctuations and reverse power flows. These issues may be mitigated with onsite energy storage systems but the latter are usually not available or expensive [2]. An alternative solution is demand side management (DSM) strategies, which may have the dual effects of reducing electricity consumption during peak hours and allowing greater efficiency and flexibility for renewable integration, namely by enabling a better match between electric supply and demand [3].

During the last decade, various models have been proposed to define household demand side management strategies for minimizing electricity cost and smoothing the load profiles. However, there are only a few studies focused on impacts of DSM on distribution system operations and renewable integration. Reference [4] proposed a framework for evaluating impacts of DSM on composite generation and transmission system reliability. Reference [5] and [6] shown the economic benefits of DSM on the agricultural and industrial sectors. References [7]–[9] introduced recent researches on the effects of DSM on various aspects of power system reliability. These studies were, however, mostly done based on the bulk power systems which did not necessarily include individual load sectors.

In this paper we propose a framework to study decentralized household demand side management in a residential distribution network which consists of multiple smart homes with schedulable electrical appliances and some rooftop photovoltaic generation units. Each smart home makes individual appliance scheduling to optimize the electric energy cost according to the day-ahead forecast of electricity prices and its willingness for convenience sacrifice. Using the developed simulation model we examine the performance of decentralized household DSM and study their impacts on the distribution network operation and renewable integration, in terms of utilization efficiency of rooftop PV generation, overall voltage deviation, real power loss, and possible reverse power flows. Compared with the work in the literature, the contributions of this work include: (1) development of a multi-household simulation framework to study decentralized DSM in a residential distribution network; (2) the proposed DSM optimization model takes into account time-varying electricity prices and rooftop PV generation available onsite; (3) every smart home that participates in the DSM program aims to reduce its electricity bill with a manageable sacrifice of convenience and comfort; (4) a comprehensive comparative study is conducted to examine the impacts of DSM on the system operation and distributed renewable integration.

The rest of this paper is organized as follows. Section II presents the system modeling for the proposed study. Section III defines the optimization model of decentralized DSM with rooftop PV at each smart household. Section IV provides numerical simulation results and analyzes the impacts of different DSM schemes on grid operation and renewable integration. Finally, section V concludes the paper and discusses future works.

## II. SYSTEM MODELLING

This section describes a system model for the proposed decentralized household DSM in a residential distribution network, which considers a single residential feeder supplying a small community with 30 residential households, as shown in Fig. 1. The time-varying load profile of each household is generated by a time-series load modeling [10] we developed using real residential demand data obtained from the open-access database (OpenEI). It is assumed that the simulated residential community has up to 16 smart homes that participate the decentralized DSM program with their schedulable electric appliances. Each smart home will optimize its appliance operation schedule to save electricity costs according to the day-ahead pricing forecast, rooftop PV generation (if available) with a controlled sacrifice of homeowner's convenience or comfort.

In order to study the impacts of decentralized DSM on renewable integration, the proposed framework also assumes that each smart household will be given an opportunity to mount a rooftop solar PV panel of a rated capacity of 6kW. For the analysis simplicity and without loss of generality, we consider all the rooftop PV systems in this small residential community have similar solar insolation and produce the same amount of electric power.

It is worth noting that each smart household has a specific set of flexible electric appliances, with different power ratings and operational limits, as shown in Table 1. More details on the household DSM appliances may be found in [10]. In order to avoid creating additional peak load period resulted from the DSM load shifting, each smart home will be given a maximum demand (MD) constraint equal to the peak load of its original load profile.

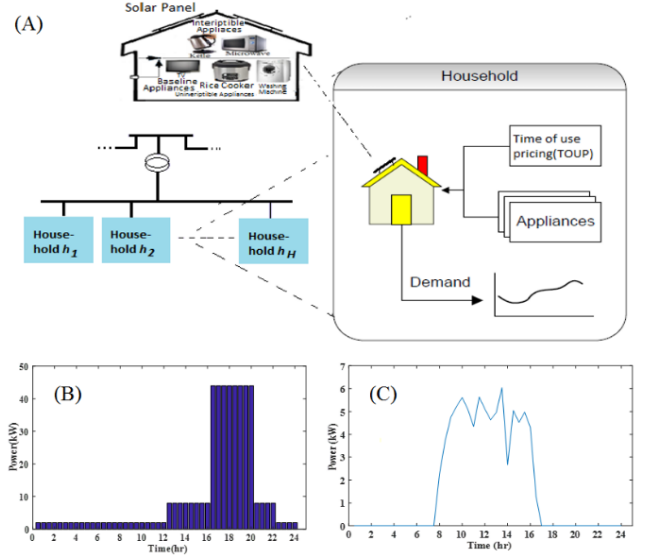

Figure 1. System model for the decentralized DSM study: (A) smart household with interruptible appliances and roof-top PV; (B) day-ahead electricity price forecast; (C) output power of rooftop PV generation.

TABLE I. SMART HOUSEHOLDS WITH AGREED MD LIMITS

| Household Index | Bus No. | Interruptible Appliance No. | Uninterruptible Appliance No. | MD(kW) |
|---|---|---|---|---|
| 1 | 12 | 21 | 7 | 12.4 |
| 2 | 14 | 15 | 4 | 15.3 |
| 3 | 17 | 17 | 4 | 11.8 |
| 4 | 20 | 12 | 4 | 8.0 |
| 5 | 3 | 21 | 7 | 12.4 |
| 6 | 10 | 15 | 4 | 15.3 |
| 7 | 19 | 17 | 4 | 11.8 |
| 8 | 29 | 12 | 4 | 8.0 |
| 9 | 8 | 21 | 7 | 12.4 |
| 10 | 15 | 15 | 4 | 15.3 |
| 11 | 22 | 17 | 4 | 11.8 |
| 12 | 25 | 12 | 4 | 8.0 |
| 13 | 6 | 21 | 7 | 12.4 |
| 14 | 16 | 15 | 4 | 15.3 |
| 15 | 23 | 17 | 4 | 11.8 |
| 16 | 31 | 12 | 4 | 8.0 |

## III. MATHEMATICAL FORMULATION

The proposed decentralized household DSM model will be implemented at each smart home participating in the program individually. It aims to minimize the household's electricity cost by scheduling the on/off status of domestic appliances over the operational periods, considering the dynamic electricity prices, locally available PV generation, and the penalty cost of appliance operation time-shifting. The penalty cost is included in order to manage the customer inconvenience caused by the DSM load shifting. Assume that the proposed demand side management program is scheduled day-ahead with 30 min per slot. The decision variables are the operational status of appliances $u_a(t)$ over the next 24 hours for each household.

Three levels of DSM participation will consider the enrollment of 4, 8, and 16 smart homes in DSM program respectively. In addition the simulation model will use three penalty prices of 0, 5, and 10 (¢/kWh) to represent different compensations requested for the sacrifice of convenience. After implementation of DSM on selected households final load profiles represent the optimal load profile of the households are considered to run AC power flow and examine the voltage, power loss and power flow across the residential feeder. The impacts of DSM on the residential network will be examined with two scenarios: (a) DSM households without PV Installation; and (b) DSM households with rooftop PV installation on site.

The decentralized DSM optimization for each household can be defined as below:

$$\min_{\{u_a(t)\}} C_e + C_p \quad (1)$$

subject to:
$$C_e = 0.5 \times \sum_{t=1}^{T}(P_{load}(t)) \times \pi_e(t) \quad (2)$$
$$C_p = 0.5 \times \pi_p \sum_{a=1}^{A} \Delta T_a \, r_a \quad (3)$$
$$P_{load}(t) = \max((\sum_{a=1}^{A} r_a \times u_a(t) - \alpha \times P_{pv}(t)), 0) \quad (4)$$
$$\sum_{a=1}^{A} r_a \times u_a(t) \leq MD \qquad \forall a \in \{1 \text{ to } A\} \quad (5)$$
$$\sum_{t=1}^{T} u_a(t) = D_a \qquad \forall a \in \{1 \text{ to } A\} \quad (6)$$
$$u_a(t) = 0 \qquad \forall t < s_a \text{ or } \forall t > f_a \quad (7)$$
$$\Delta T_a = 1^T \cdot |t_a^{st_{new}} - t_a^{st_{old}}| \qquad \forall a \in \{1 \text{ to } A\} \quad (8)$$
$$t_a^{st_{new}} = [t|u_a^{new}(t) = 1]_{1 \times D_a} \qquad \forall a \in \{1 \text{ to } A\} \quad (9)$$
$$t_a^{st_{old}} = [t|u_a^{old}(t) = 1]_{1 \times D_a} \qquad \forall a \in \{1 \text{ to } A\} \quad (10)$$

where eq. (2) and (3) define the electricity cost and the penalty cost of DSM scheduling respectively. The proposed model assumes that surplus PV generation will be injected into the distribution network without reward. Therefore, the total electricity cost within each time slot should be no less than zero as indicated by eq. (4). The binary factor α in eq. (4) represents the status of PV installation at a DSM household. Constraint (5) presents the Maximum Demand (MD) that the aggregate appliance power of each household cannot exceed at any time. This specified upper limit is to prevent super-high power demand peak even during the hours when day-ahead electricity price is low because the utilities do not want to have "new" peak created by the DSM load-shifting or because the distribution feeders have capacity constraints. Constraint (6) and (7) indicate the total operation duration and the allowable turn-on time of an appliance and constraints (8)-(10) specify the original and the new starting point, $t_a^{st_{old}}$ and $t_a^{st_{new}}$ respectively, to capture the duration of time-shifting for flexible appliances.

## IV. NUMERICAL SIMULATION RESULTS

Numerical experiments have been conducted to examine the impacts of decentralized household DSM program with rooftop PV system on the residential distribution network in terms of renewable usage efficiency, voltage deviations, real power loss, and reverse power flow. In the first scenario the proposed DSM model is considered without rooftop PV installation on selected DSM households and the impact of different DSM penetration level with different penalty prices are investigated. In second scenario, DSM households considered to have rooftop PV systems and with calculation of optimal load profiles for DSM households, distribution network operation in terms of voltage, power loss and reverse power flow is investigated.

### A. Individual Household DSM

Fig. 2 compares the original load profiles (depicted as red solid line) and the DSM scheduled profile (depicted as blue bar plots) of smart home #1 under four different conditions (with/without onsite PV, and $\pi_p = 0$ or 5 ¢/kWh). The PV output power is depicted as blue dashed line. Original load profile indicates two demand peaks: around 7:00 to 9:00 in the morning and 17:00 to 20:00 in the evening.

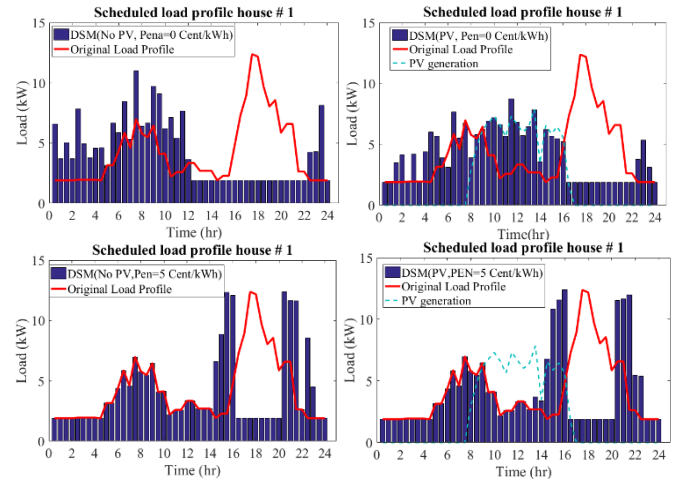

Figure 2. Load profiles of smart household #1 with DSM under four different conditions

As shown in Fig.1, there is a peak price rate during the time of 17:00 to 20:00. Therefore, the DSM program without any penalty for load-shifting inconvenience will move all the schedulable appliances out of the peak price hours, as shown in the top two subfigures in Fig. 2. When the onsite rooftop PV is mounted, the appliance load will be first shifted to the time slots inside 8:00-16:00 where the PV generation is available hence achieve a renewable usage efficiency of 99.98%. However, if penalty for load-shifting inconvenience is considered, say, with $\pi_p = 5$ ¢/kWh, the shifted load will be concentrated on the boundary next to peak-price hours as 14:00-16:00, and 20:00-22:00. And this may significantly affect the renewable usage efficiency, causing it to drop to 36.55%.

### B. Voltage Fluctuations

Fig. 3 compares the voltage profile at the end of the feeder for different DSM participation levels with $\pi_p = 0$ ¢/kWh considering the two scenarios, i.e., with or without PV installation. It is found that higher DSM participations will tend to "flatten" the voltage profile (i.e. filling the voltage drop valleys) and bring better improvement to voltage fluctuations caused by load variations. Besides, during the daytime hours (8:00-16:00) when the rooftop PVs generate power, the DSM without penalty for convenience sacrifice will help mitigate overvoltage problem during those hours.

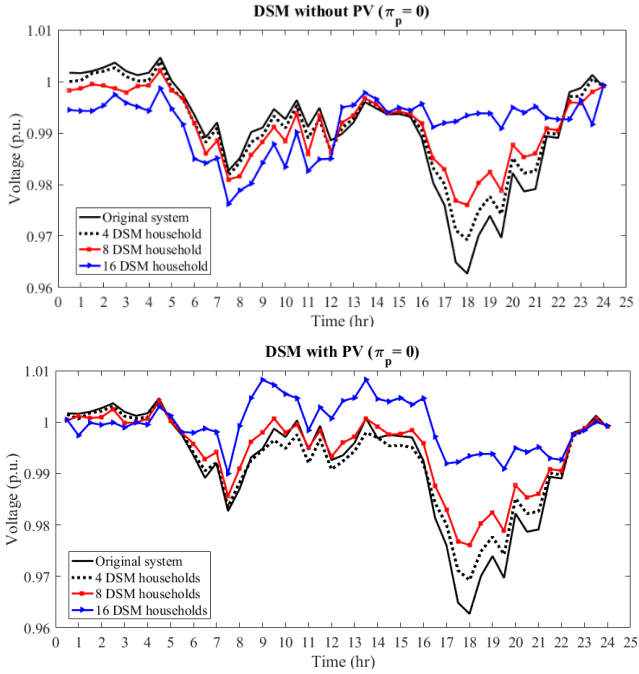

Figure 3. Voltage profile at the end of the feeder for different DSM participation levels (with/without PV installation, $\pi_p = 0$).

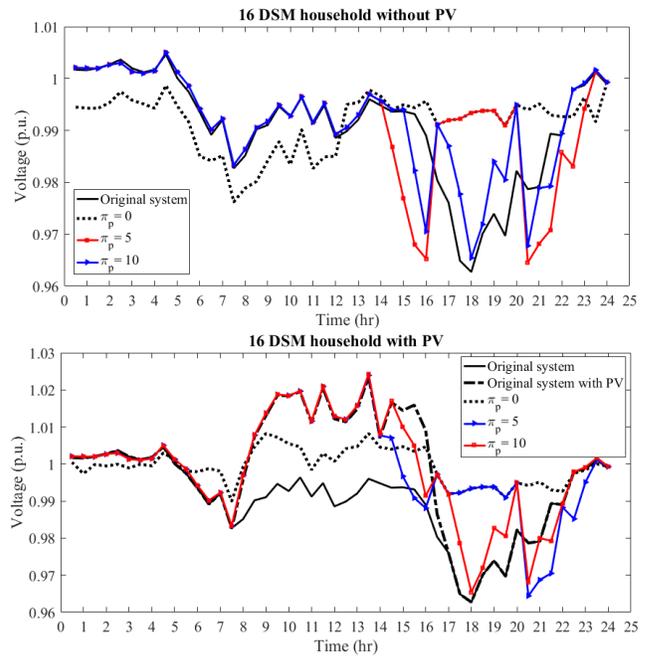

Figure 4. Voltage profile at the end of the feeder with 16 DSM households with different penalty prices

Fig.4 compares the voltage profile at the end of the feeder when there are 16 smart homes participating in the DSM program with different penalty prices considering the two scenarios, i.e., with or without PV installation. Comparing with the results in Fig.3, we may see that higher penalty prices such as $\pi_p = 5$ or 10 ¢/kWh will significantly reduce the voltage smoothing effect of the DSM program because much fewer flexible appliances will be shifted out from the peak-price hours. Besides, during the daytime hours (8:00-16:00) when the rooftop PVs generate power and cause overvoltage conditions, the DSM with high penalty for convenience sacrifice cannot help much during those hours.

## C. Reverse Power Flow

Fig. 5 presents the power flow distribution across the distribution feeder of the original system and with 16 household DSM ($\pi_p = 0$), without any PV installation. Obviously, DSM participations without high penalty for load shifting inconvenience help smoothing the flow distribution. Fig. 6 shows the power flow distribution across the distribution feeder of the original system and with 16 household DSM ($\pi_p = 0$ ¢/kWh ), with 16 onsite rooftop PV installations. The output power from rooftop PVs may cause reverse flows at some line segments at mid-day hours as can be seen from the top subfigure. However, the DSM scheduling may shift some appliance usage into those time slots and lessen the reverse flows.

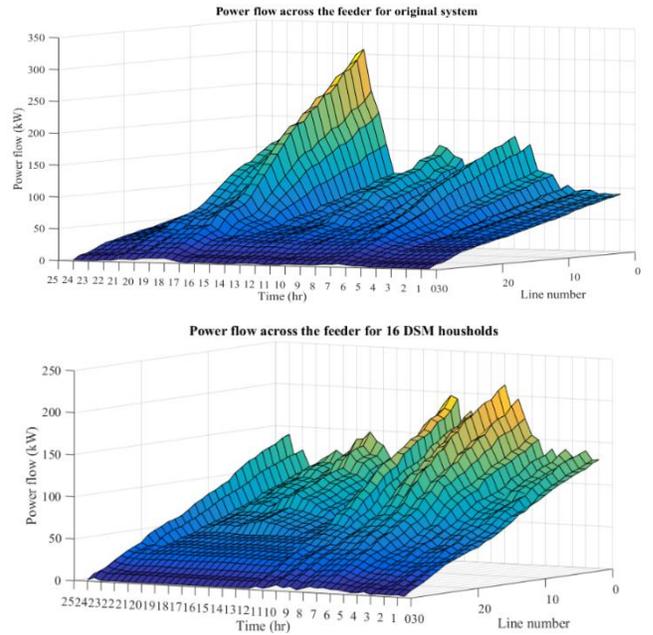

Figure 5. Power flow across the feeder without PV installation

## D. Real Power Loss

Fig. 7 compares the feeder real power loss experienced at each hour with different DSM participation levels with or without PV installation. Clearly, it shows that higher DSM participations will more effectively bring down the high energy loss during peak-price hours with a little bit increase during the mid-day hours. However, with onsite rooftop PV installations, the output PV generated power during the mid-day time will help to get rid of the real-power loss caused by DSM load shifting.

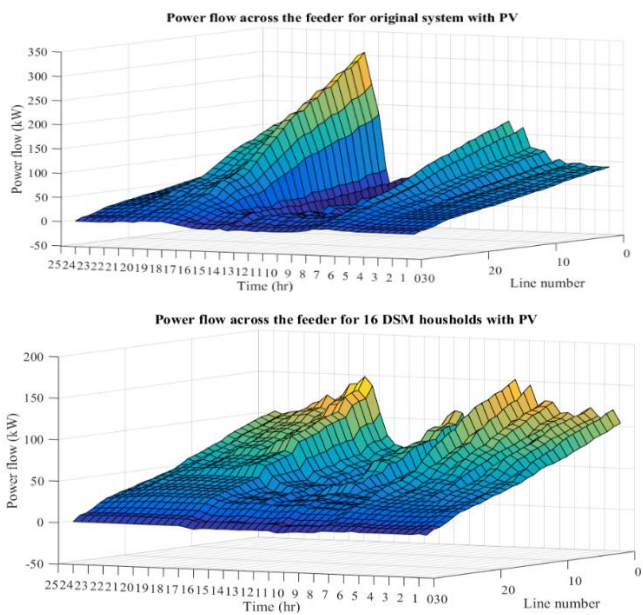

Figure 6. Power flow across the feeder with on-site rooftop PV installations

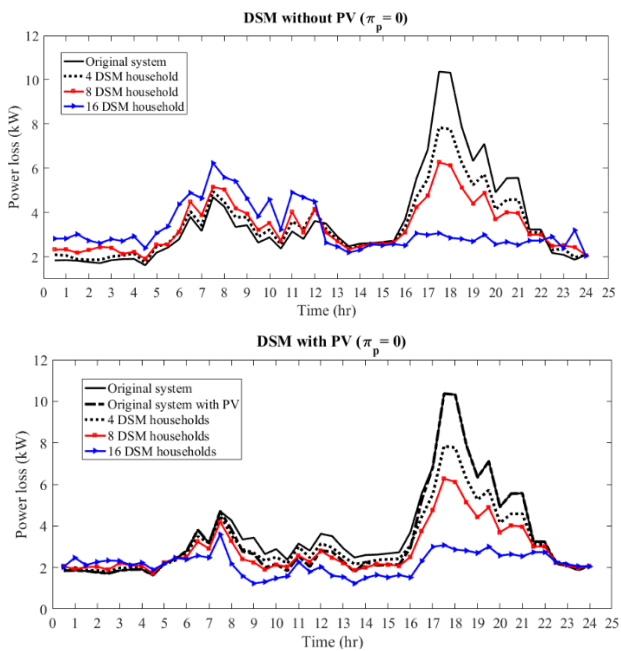

Figure 7. Power loss in the feeder for different DSM participation levels with or without PV installation

## V. CONCLUSIONS AND FUTURE WORK

This paper has developed an efficient DSM simulation framework for a residential distribution network with multiple smart homes and onsite PV installation. Numerical experiments are performed to examine the interactions of dynamic electricity price signal, the penalty for load shifting inconvenience, the load shifting of flexible appliances in a household, the PV generation power from rooftop PV installations. The performance of decentralized household DSM and their impacts on the distribution network operation and renewable integration haven been studied in terms of utilization efficiency of rooftop PV generation, overall voltage deviation, real power loss, and possible reverse power flows. It is found that DSM without high penalty for load-shifting inconvenience will most effectively smooth the load profile, reduce the voltage fluctuations, encourage renewable generation consumption, and mitigate the overvoltage and reverse power flow problems caused by PV generated power during the mid-day. In the future, we plan to augment the DSM simulation framework so it may be applied for both centralized and decentralized DSM program for multi-category customers such as residential, commercial, and/or industrial. We are also interested to study the fast dynamics of renewable distributed generation and their impacts on the voltage regulations and transient stability of distribution network.

REFERENCES

[1] International Energy Agency (2011). Technology Roadmap. Smart Grids. [Online] Available: http://www.iea.org/publications/frepublications/publication_roadmap.pdf.
[2] A. Shah, Energy security, accessed on Nov.2016 [Online]. Available: http://www.globalissues. Org/article/595/energy-security.
[3] IEA. *Strategic plan for the IEA Demand –Side Management program 2008-2012*. Accessed on Nov. 2016. [Online]. Available: http://www.ieadsm.org.
[4] S.H. Elyas, H. Sadeghian, Hayder O. Alwan, Zhifang Wang "Optimized Household Demand Management with Local Solar PV Generation" 2017 North America Symposium(NAPS),Sep 17-19, 2017.
[5] M. Motalleb, A. Eshraghi, E. Reihani, H. Sangrody, and R. Ghorbani, "A Game-Theoretic Demand Response Market with Networked Competition Model," in 49th North American Power Symposium (NAPS), 2017.
[6] C. Goldman, M. Reid, R. Levy, and A. Silverstein, "Coordination of energy efficiency and demand response," Ernest Orlando Lawrence Berkely Nat. Lab., Jan 2010. [Online]. Available: https://emp.lbl.gov.
[7] P. Cappers, C. A. Goldman, and D. Kathan, "Demand response in U.S. electricity markets: Empirical evidence," *Energy, vol. 35, no. 4, pp. 1526-1535, 2010.*
[8] J. h. Yoon, R. Bladick, and A. Novoselac, "Demand Response for residential buildings based on dynamic price of electricity," *Energy Build.*, vol. 80 pp, 531-541, Sep. 2014.
[9] M. Paulus and F. Borggrefe, "The potential of demand-side management in energy intensive industries for electricity markets in Germany,"Appl. Energy, vol. 88, pp. 432–441, 2011.
[10] H. Sadeghian, M. H. Athari, Zhifang Wang "Optimized Solar Photovoltaic Generation in a Real Local Distribution Network". *IEEE PES 2017 Innovative Smart Grid Technologies (ISGT 2017) Conference, Arlington,* VA USA, Apr. 23-26, 2017.